\documentclass{article}
\usepackage[leqno]{amsmath}
\usepackage{amsfonts}
\usepackage{amsthm}
\usepackage{amsmath}
\usepackage{amssymb}
\usepackage{amscd}
\usepackage{inputenc}
\usepackage{t1enc}
\usepackage[mathscr]{eucal}
\usepackage{indentfirst}
\usepackage{graphicx}
\usepackage{graphics}
\numberwithin{equation}{section}
\usepackage{epstopdf} 
\usepackage{multirow}
\usepackage{mathtools}
\usepackage{caption}
\usepackage{setspace}
\usepackage{tikz-cd}
\usepackage{multicol}
\usepackage{comment}
\usepackage[title,titletoc]{appendix}

\usepackage[style=alphabetic, backend=biber]{biblatex}
\usepackage[cache=false]{minted}
 
\newcommand{\bibhack}[1]{}

\theoremstyle{plain}
\newtheorem*{theorem*}{Theorem}
\newtheorem*{cor*}{Corollary}
\newtheorem*{lem*}{Lemma}
\newtheorem{Th}{Theorem}[section]
\newtheorem{Lemma}[Th]{Lemma}
\newtheorem{lemma}[Th]{Lemma}

\newtheorem{cor}[Th]{Corollary}

\newtheorem{maintheorem}{Theorem}
\newtheorem{maincor}[maintheorem]{Corollary}

 \theoremstyle{definition}
\newtheorem{Def}[Th]{Definition}
\newtheorem*{Def*}{Definition}

\newtheorem{remark}[Th]{Remark}

\newtheorem{?}[Th]{Problem}

\theoremstyle{plain}
\newtheorem{innercustomgeneric}{\customgenericname}
\providecommand{\customgenericname}{}
\newcommand{\newcustomtheorem}[2]{%
  \newenvironment{#1}[1]
  {%
   \renewcommand\customgenericname{#2}%
   \renewcommand\theinnercustomgeneric{##1}%
   \innercustomgeneric
  }
  {\endinnercustomgeneric}
}

\newcustomtheorem{customthm}{Theorem}
\newcustomtheorem{customtheorem}{Theorem}
\newcustomtheorem{customlemma}{Lemma}
\newcustomtheorem{customcor}{Corollary}
\newcustomtheorem{customconj}{Conjecture}

\newcommand{\Nu}{\mathcal{V}}

\newcommand{\setvert}{\,\vert \,}

\newcommand{\Q}{\mathbb{Q}}
\newcommand{\Z}{\mathbb{Z}}

\newcommand{\C}{\mathbb{C}}

\newcommand{\cp}{\mathbb{CP}}
\newcommand{\CP}[1]{\cp^{\frac{#1}{2}}}
\newcommand{\hp}{\mathbb{HP}}
\newcommand{\HP}[1]{\hp^{\frac{#1}{4}}}
\newcommand{\into}{\xhookrightarrow{}}
\newcommand{\acton}{\curvearrowright}

\newcommand{\ft}{\mathfrak{t}}
\newcommand{\PP}{\mathbb{P}}
\newcommand{\OO}{\mathcal{O}}
\newcommand{\Fr}{{\rm Fr}}
\newcommand{\norm}[1]{\left\lVert#1\right\rVert}
\renewcommand{\H}{\mathcal H}
\newcommand{\bb}{\mathbb}
\newcommand{\Tr}{{\rm Tr}}
\newcommand{\V}{\mathcal V}
\newcommand{\orb}{{\rm orb}}
\renewcommand{\O}{{\rm O}}
\newcommand{\orient}{{\rm or}}

\newcommand{\tensor}{\otimes}

\newcommand{\lift}[1]{\overline{#1}}
\newcommand{\res}[1]{#1|_C}
\newcommand{\lr}[1]{\res{\lift{#1}}}

\title{The Hopf conjecture with $T^3$-symmetry}
\author{Jan Nienhaus, Dennis Wulle}
\date{September 2026}

\begin{document}

\maketitle
\begin{abstract}
    We prove that any closed even-dimensional manifold with positive sectional curvature admitting an isometric action by a 3-torus has positive Euler characteristic.  We also prove the result for orbifolds.
\end{abstract}

\section{Introduction}
In the 1930s, Hopf conjectured that every even-dimensional manifold with positive sectional curvature has positive Euler characteristic. While little pro\-gress has been made on the conjecture in full generality, as part of the \textit{Grove Symmetry Program}, initiated in the 1990s, significant advancements have occurred in showing the conjecture holds true with ever-smaller symmetry assumptions. This may be even more surprising considering there are some recent signs that the conjecture may be false in full generality\footnote{There has been an announcement of a counterexample to the Hopf conjecture in dimension 6, but to the knowledge of the authors this work is entirely AI and not verified by geometers. The Hopf conjecture is known to hold for $6$-manifolds admitting an isometric circle action \cite{PuttmannSearle2001}.}. We prove

\begin{maintheorem}\label{main:Hopf}
    Let $M$ be a closed positively curved manifold of even dimension, admitting an effective isometric $T^3$-action. Then $\chi(M)>0.$
\end{maintheorem}

The Euler characteristic of $M$ is known to be equal to that of the necessarily nonempty fixed point set of the action. Thus Theorem \ref{main:Hopf}  follows from

\begin{maintheorem}\label{main:Fixedpointclassification}
    Let $M$ be an oriented closed positively curved manifold, admitting an effective isometric $T^3$-action. Let $F$ be any fixed point component of such an action, with dimension $f$. Then the rational cohomology ring of $F$ is that of $S^f, \CP{f}, \HP{f},$ or $S^2\times\HP{f-2}$. 
\end{maintheorem}

This represents a capstone for what one might call the $4$-periodicity program for the Hopf conjecture. A cohomology ring is called $4$-periodic if it's spherical or there is a class in degree $4$ such that cupping with that class gives isomorphisms in all degrees. All the model spaces of Theorem \ref{main:Fixedpointclassification} are examples of this; in fact establishing $4$-periodicity of $M$ or its fixed point components has been part of all major work on the problem in the last 25 years.

In terms of rank, Theorem \ref{main:Fixedpointclassification} is sharp. The positively curved $24$-dimensional octonionic flag manifold contains a $6$-dimensional complex flag manifold as $T^2$-fixed point component. This does not have one of the rational cohomology rings of Theorem \ref{main:Fixedpointclassification}, in fact it has Betti numbers $b_2=b_4=2$, while all models in Theorem \ref{main:Fixedpointclassification} have $b_i\le 1$ for all $i$. In particular $b_4\ne b_0=1$, so this ring is not $4$-periodic in the strict sense.

Periodicity in this context was introduced in \cite{Wil03}. Under a linear lower bound $T^d\acton M^n$, $d\gtrsim\frac{n}{6}$, using his new Connectedness Lemma, Wilking is able to show the cohomology of $M$ itself is one of a list of $4$-periodic models. 

The next breakthrough in the theory was \cite{Kennard13}. Using mod-$p$ methods, Kennard is able to prove automatic upgrading of $\Z$-periodicity with any length to $4$-periodicity over $\Q$, which allows him to prove $4$-periodicity for fixed point components assuming only $d\gtrsim \log(n)$. To use his theorem, one needs to produce configurations of submanifolds satisfying certain codimension bounds. 

While his original work proves this is possible by induction, Kennard, Wie\-me\-ler and Wilking \cite{Kennard2021-ge} were later able to show, using a newly developed splitting theorem for torus representations, that the required situation can be arranged without using induction, which results in a cohomology classification for fixed point components assuming only a $T^5$-action. 

By proving a $4$-periodicity theorem under much sharper assumptions, the first named author was later able to prove the same classification using only a $T^4$-action \cite{Nienhaus22}. The significant difference is that the generic mod-$p$ methods of Kennard are replaced by an argument passing entirely through characteristic classes. 

This characteristic class approach also drives our proof of Theorem \ref{main:Fixedpointclassification}. However, we will be able to directly show that a specific rational degree $4$ class (the first Pontryagin class of a suitable vector bundle) induces $4$-periodicity on $F$, completely bypassing the previous mod-$p$-methodology and $4$-periodicity theorems.

Under the action of a $4$ resp. $5$-torus, the cohomology ring of $\bb S^2\times \HP{f-2}$ is ruled out in \cite{Nienhaus22} and \cite{Kennard2021-ge}. It remains an interesting question whether this case can occur for $3$-torus actions or not. All other cases are known to occur. 

Since there are vastly more known examples of non-smooth spaces of positive curvature compared with smooth ones, one may ask whether the conclusions of Theorems \ref{main:Hopf} and \ref{main:Fixedpointclassification} depend crucially on smoothness. The natural non-smooth test case is orbifolds (for Alexandrov spaces, every imaginable version of Theorem \ref{main:Fixedpointclassification} is wrong because one has arbitrary suspensions).

In fact, by establishing a Connectedness Lemma for orbifolds, were are able to generalize Theorems \ref{main:Hopf}
 and \ref{main:Fixedpointclassification} to the non-smooth case:

 \begin{maintheorem}
     Let $\OO$ be an orientend closed positively curved orbifold admitting an effective isometric $3$-torus action. Let $F$ be any fixed point component of such an action with dimension $f$. Then the rational cohomology ring of the underlying topological space $\vert F \vert$ of $F$ is that of $\bb S^f$, $\CP{f}$, $\HP{f}$ or $\bb S^2 \times \HP{f-2}$.
 \end{maintheorem}

 \begin{maincor}\label{main:orbiHopf}
     Let $\OO$ be a closed positively curved orbifold of even dimension admitting an effective $3$-torus action. Then $\chi(\vert \OO \vert) >0$.
 \end{maincor}

Orbifolds with positive curvature have seen increasing amounts of consideration in recent years, see \cite{Florit2009, Gorodski2016, WulleZarei26, NienhausWulle2026}. Besides taking the Euler charactersistic of the underlying space, there are two other common notions of Euler characteristic for orbifolds, the Euler-Satake characteristic introduced in \cite{Satake1957} and the stringy Euler characteristic \cite{ALR07}. It would be interesting to know whether Corollary \ref{main:orbiHopf} holds for the other notions too. This is true for the known examples.

A few words on the proof of Theorem \ref{main:Fixedpointclassification}: Fixing a fixed point component $F$, the first step is to use the Sylvester-Gallai theorem on the weights of the isotropy action at $F$ in order to find a circle $S^1\subset T^3$ such that the remaining action $T^2=T^3/S^1 \acton Fix(S^1)\subset M$ splits near $F$ as a product of circle actions, i.e. such that at a point in $F$ the weights of the isotropy representation lie on two lines. This is in the spirit of the main reduction in \cite{Kennard2021-ge}, noting that rank three is low enough that no general machinery is necessary. Now the main technical result of the paper is

\begin{maintheorem}\label{main:t2}
    Let $M'$ be closed, positively curved, and oriented, with an effective isometric action of $T^2$ and assume $F$ is a fixed point component such that the weights of the isotropy action of $T^2$ at $F$ lie on the union of two lines. Then $H^*(F,\Q)$ is $4$-periodic with vanishing odd cohomology in degrees $0<*<\dim(F)$. This leaves only the rings in Theorem \ref{main:Fixedpointclassification}.
\end{maintheorem}

The two lines of weights correspond to two circles $S^1_\alpha, S^1_\beta\subset T^2$ whose fixed point components containing $F$ we will call $F_\alpha, F_\beta$. The condition on the weights make $F_{\alpha, \beta}$ intersect transversely with intersection $F$. By Wilking's connectedness Lemma \cite{Wil03}, the inclusion $F\into F_\alpha$ is $\dim(F)$-connected. Our argument passes through Borel-equivariant cohomology $H^*_T$. A standard spectral sequence argument shows that the connectedness passes over to give that $H^*_T(F_\alpha)\to H^*_T(F)$ is an isomorphism in degrees $*<\dim(F)$ and injective in degree $\dim(F)$. In $H^{\dim(F)}_T(F)$ there is a class $[F]$ corresponding to the orientation class of $F$. If $F$ has cohomology ring different from one of the examples, $[F]$ lifts to $H^*_T(F_\alpha)$ and the Atiyah-Bott localization formula shows that there is a second fixed point component $C$ such that the restriction of this lift to $C$ is a non-vanishing zero divisor in $H^*_T(C)$.

Finally, one considers the first Pontryagin class of the normal bundle $\Nu$ of $F_\alpha$, restricted to $F$. One can show that if cup product with this class on $F$ is not injective in every degree, the restriction of this class to $C$ is also a zero divisor. A class in $H^*_T(C)=H^*(C)\otimes H^*(BT)$ is a zero divisor exactly if it has no part in $H^0(C)\otimes H^*(BT)$. Vanishing of this component for $p_1(\Nu|_F)$ at $C$ then gives a relation between the weights of the isotropy action at $F$ and $C$, which forces the same structure on these weights at $C$ as they have at $F$. This finally results in a contradiction to Frankel's theorem, as it forces $F$ and $C$ to intersect while they were constructed to be disjoint.

\subsection*{Acknowledgements}
The first named author would like to thank the Cluster of Mathematics M\"unster for enabling a research stay in M\"unster in the summer of 2026, where part of this work was carried out.

The second named author acknowledges support by the Alexander von Humboldt Foundation through Gustav Holzegel’s Alexander von Humboldt Professorship endowed by the Federal Ministry of Education and Research and by the Deutsche Forschungsgemeinschaft (DFG, German Research Foundation) under Germany's Excellence Strategy EXC 2044/2 –390685587, Mathematics Münster: Dynamics–Geometry–Structure.

\subsection*{AI disclosure}
The starting point for this paper was a proof of Theorem \ref{main:Hopf} obtained with OpenAI's ChatGPT-6-Astra and minimal human input, using the $c_1$-approach outlined in Appendix A. This proves vanishing of the odd cohomology of the fixed point components, but no additional structure. The alternative $p_1$-route was discovered interactively when trying to eliminate the finite-isotropy reduction from the argument. 
The observations that the new route allows for structure results on the cohomology of fixed point components and the generalization to orbifolds are fully human. In the orbifold part AI tools have been used for literature research.  

The structuring of the argument, including all statements and proofs as appearing in this article, was done by the authors. This document contains no machine-generated text.

\section{Manifolds}

We proceed by an important reduction, following the approach taken in \cite{Kennard2021-ge}:\\[2mm]
\noindent
\textbf{Reduction $\Delta$:}
Instead of an action $T^3\acton M$, we may consider an action $T^2=T^3/S^1\acton M'=M^{S^1}$ such that at a point $p\in F$ the isotropy action only has weights that are multiples of two linearly independent weights $\alpha, \beta$.

\begin{proof}
    Take the set of weights $\Omega\subset\ft^*$. Projectivize and use Sylvester-Gallai to obtain a line in $\PP^2=\PP(\ft^3)$ that contains exactly two points of $\PP(\Omega)$. 
    This gives a $2$-dimensional subspace $V$ of $\ft^3$ such that $\Omega\cap V$ is contained on the union of two distinct lines through the origin. 
    Let $S^1$ be the subcircle of $T^3$ with Lie algebra $\ker(V)$. Then $(T^2=T^3/S^1\acton M^{S^1}, F)$ is as claimed.
\end{proof}

\begin{remark}
    Replacing the $S^1$  by a possibly disconnected 1-dimensional subgroup, this is a special case of \cite[Theorem F]{Kennard2021-ge}. In the manifold case this would also be sufficient, but we want to avoid disconnectedness issues in anticipation for the generalizations in section \ref{sec:orbifolds}.
\end{remark}

\begin{lemma}
    Let $(T^2\acton M, F), p, \alpha, \beta$ be as in reduction $\Delta$. Let $S^1_\alpha$ be the circle with Lie algebra $\ker(\alpha)$, and $F_\alpha$ the fixed point component of $S^1_\alpha$ containing $F$; same for $\alpha$ replaced with $\beta$. Let $a,b$ be the total multiplicity of weights that are nonzero multiples of $\alpha, \beta$ as representations in the normal bundle of $F$ at $p$.

    Then $\dim(F_\alpha)=\dim(F)+2a$, $\dim(F_\beta)=\dim(F)+2b$ and, if ordered such that $a\le b$, $F=F_\alpha\cap F_\beta$ and the inclusion $F\into F_\alpha$ is $\dim(F)$-connected.
\end{lemma}

\begin{proof}
    At any fixed point of $T$, the weights of the induced action of $S^1_\alpha$ are given by taking the weights of $T^2$ mod $\alpha$. Thus the space of fixed directions at $p$ is exactly the sum of the fixed space of $T^2$ and the irreducible summands whose weight is linearly dependent with $\alpha$.
    $F_\alpha$ and $F_\beta$ intersect transversely in $M$ by counting dimension. Thus by the Connectedness Lemma, $F_\alpha\cap F_\beta \into F_\alpha$ is $\dim(F)$-connected if $F_\alpha$ has the smaller dimension of $F_\alpha, F_\beta$. in particular, $F_\alpha\cap F_\beta$ is connected, so must equal $F$.
\end{proof}

\subsection{Preliminaries: bundle preparation}

\begin{lemma}
    Consider the normal bundle $\Nu$ of $F_\alpha$ in $M$. Then:

    $F_\alpha$ is $T^2$-invariant, thus $\Nu$ carries an induced $T^2$-action. At a point $p\in F$ the restriction $\Nu|_p$ is a direct sum of $b$ copies of representations with weights multiples of $\beta$.
\end{lemma}

\begin{proof}
    $T^2$ is abelian, thus fixes the fixed point components of any subgroup, in particular $F_\alpha$. The statement at $p$ is exactly the claim of reduction $\Delta$, noting that the irreducible subrepresentations with weight $c\alpha$ lie in the tangent bundle of $F_\alpha$ and those with weight $c\beta$ are normal to $F_\alpha$.
\end{proof}

\begin{lemma}
    The bundle $\Nu$ over $F_\alpha$ can be equipped with a complex structure compatible with the $T$-action such that all $\beta$-type subrepresentations at $p$ come with the same sign.
\end{lemma}

\begin{proof}
    To fix a complex structure on $\Nu$ it is enough to fix an oriented circle in $T^2$ that fixes the base point of $\Nu$ pointwise and acts without kernel on the fiber: Given such a circle and a base point $x$, split $\Nu_x$ by weights, and on each part declare that $J(v)=\frac{\xi v}{|\xi v|}$ where $\xi\in\ft$ is the positively oriented generator of the Lie algebra of the chosen $S^1$. This is easily checked to be a complex structure, and the construction commutes with the action of $T$, so the action of $T$ on $\Nu$ becomes $\C$-linear.

    Now do this with $S^1=S^1_\alpha$. By construction this fixes $F_\alpha$ pointwise and does not fix any direction normal to $F_\alpha$, so the construction applies. Furthermore, the definition is such that $S^1_\alpha$ always acts in positive orientation with respect to $J$, so all occurrences of $\beta$ at $p$ project in $\ft$ to positive multiples of the chosen generator of the Lie algebra of $S^1_\alpha$, in particular must be positive multiples of each other.
\end{proof}

\subsection{Preliminaries: equivariant cohomology}\label{prem:eqcohom}

\begin{lemma}\label{lem:equivconnectedness}
    For $T\acton X$, denote by $X_T$ the Borel construction. Then the inclusion $F_T\into ({F_\alpha})_T$ is $\dim(F)$-connected.
\end{lemma}

\begin{proof}
    This is the 5-Lemma for the induced map of homotopy exact sequences for the fibrations $(F\to F_T\to BT)\into (F_\alpha \to (F_\alpha)_T\to BT)$.
\end{proof}

\begin{lemma}\label{p1lemma}
    Let $V\to F$ be a $T$-equivariant complex vector bundle with the $T$-action fixing $F$ pointwise. For each weight $\chi$ of $T^2$, denote by $V_\chi$ the subbundle of $V$ where $T$ acts on each fiber via $\chi$.
    Then the equivariant first Pontryagin class satisfies $$p_1^T(V)=p_1(V)\otimes 1 + 2\sum c_1(V_\chi) \otimes \chi + \sum \dim_\C(V_\chi) 1\otimes \chi^2$$ in $H^*(F)\otimes H^*(BT)=H_T^*(F)$, where $c_1$ is the non-equivariant first Chern class, $p_1$ is the non-equivariant first Pontryagin class, and the sums range over the weights of all irreducible subrepresentations of $V$ at any point $p\in F$, with the lattice of weights identified with $H^2(BT, \Z)\subset H^2(BT, \Q)$.
\end{lemma}

\begin{proof}
    This follows from the splitting principle. If $V_\chi$ splits as a sum of line bundles $\bigoplus L_i$ with first Chern classes $y_i$, then one has $c_1=\sum y_i$ and $p_1=\sum y_i^2$. The result is then a direct consequence of the fact that the equivariant first Chern class of $L_i$ is $y_i+\chi$, so $p_1^T=\sum(y_i+\chi)^2=\sum y_i^2+2(\sum y_i)\chi+\sum\chi^2$, together with the fact that $p_1^T(\bigoplus V_\chi)=\sum p_1^T(V_\chi)$, which itself follows directly from the formula for sums of line bundles.
\end{proof}

\begin{remark}
    At first glance, the above formula looks like it depends on the choice of complex structure. This is in fact not the case: Any other complex structure compatible with the $T$-action will flip the orientation of some subset of the $\chi$, which will give both the $\chi$ and the corresponding $c_1$ a factor of $-1$, leading to the same overall cohomology class.
\end{remark}

\subsection{Proof of the theorem}

Assume we are in the situation after reduction $\Delta$. We proceed by a sequence of lemmas, but first, some notation:

\begin{Def}
    For a class $x\in H^*(F)$, write $\lift{x}$ for the preimage of $x\otimes 1\in H^*_T(F)=H^*(F)\otimes H^*(BT)$ under the map $H^*_T(F_\alpha)\to H^*_T(F)$.

    For a class $y$ in $H^*_T(F_\alpha)$ and a fixed point component $C$, we write $\res{y}$ for its equivariant restriction in $H^*_T(C)$.
\end{Def}

By Lemma \ref{lem:equivconnectedness}, a preimage $\lift{x}$ exists and is unique if $\deg(x)\ne\dim(F)$, and is unique if it exists if $\deg(x)=\dim(F)$.

\begin{lemma}\label{orientationlifting}
    If $F$ is not a rational homology sphere, the orientation class $[F]=[F]\otimes 1$ lies in the image of $H^*_T(F_\alpha)\into H^*_T(F)=H^*(F)\otimes H^*(BT)$. Furthermore there is at least one other fixed point component $C\subset F_\alpha$ where the restriction $\lr{[F]}$ of $\lift{[F]}$ is nonzero as an element of $H^*_T(C)$.
\end{lemma}

\begin{proof}
    If $F$ is not a homology sphere, we can find some nonzero class $x\in H^*(F)$ with $0<*<\dim(F)$. Let $y$ be Poincar\'e-dual to $x$. Then $[F]=x \cup y$. Both $x$ and $y$ live in degrees where $H_T^*(F_\alpha)\to H_T^*(F)$ is an isomorphism, so we may pick lifts $\lift{x}, \lift{y}$ of $x\otimes 1, y\otimes 1$, and see that $\lift{x}\cup \lift{y}$ is a lift of $[F]$.

    For the second claim, consider the "integrate-out-the-fiber"-map $$\int^T_{F_\alpha}:H^*_T(F_\alpha)\to H^{*-\dim(F_\alpha)}_T(point)=H^{*-\dim(F_\alpha)}(BT).$$ Evaluated on the class $\lift{[F]}$ we obtain that $\int^T_{F_\alpha}\lift{[F]}\in H^{f-(f+2a)}(BT)=0$ vanishes because the group is zero. On the other hand, the Atiyah-Bott localization theorem says that after tensoring as $H^*(BT)$-modules with the fraction field $k=Frac(H^*(BT))$, there is a formula \begin{align*}
        \int^T_{F_\alpha}u = \sum\limits_{C\subset F_\alpha}\int^T_C\frac{u|_C}{e^T(\nu(C\subset F_\alpha))}
    \end{align*}
    where the sum ranges over the fixed point components of $T$ in $F_\alpha$, and the denominators are units in the localized ring. \cite[eq. 3.8, preceding discussion for unit claim]{ATIYAH19841}.

    Computing the term in this formula for $C=F, u=\lift{[F]}$ gives $$\int^T_F\frac{\lift{[F]}|_{F}}{e}=\frac{1}{e}\ne 0,$$

    where $e=e^T(\nu(F\subset F_\alpha))$ and the $1$ in the last step is $1\in H^0(BT)$. In particular some other term in the localization formula must give nontrivial contribution, so it is impossible that $\lr{[F]}=0$ for all $C\ne F$, which was the claim.
\end{proof}

\begin{lemma}\label{Fzerodivisor}
    Let $F$ and $C$ be as in Lemma \ref{orientationlifting}. If $F$ has $b_i\ne 0$ for any $0<i<f$ odd, or if there are nonzero classes $x,y\in H^*(F)$ with $\deg(x)+\deg(y)\le \dim(F)$ and $xy=0$, then $\lr{[F]}$ is a zero divisor in $H^*_T(C)$. In the second case, so are $\lr{x}, \lr{y}$.
\end{lemma}

\begin{proof}
    We begin with the odd case. Here we can write $[F]=x\cup y$ with $x\in H^{odd}(F)$. Choose lifts $\lift{x}, \lift{y}$ of $x\otimes 1, y\otimes 1$ as before, and get $\lift{[F]}=\lift{x}\cup \lift{y}$. Restricting to $C$ gives $\lr{[F]}=\lr{x}\cup \lr{y}$. Since $x$ has odd degree, $\lr{x}\in H^{odd}(C)\otimes H^*(BT)$ cannot have parts in $H^0(C)\otimes H^*(BT)$. In particular, $\lr{[F]}\in H^{\ge 3}(C)\otimes H^*(BT)$ cannot either, noting that $H^{odd}(C)=H^{odd\ge3}(C)$ since positive curvature implies $H^1(C, \Q)=0$. Since the zero divisors of $H^*_T(C)$ are exactly the elements with no part in $H^0(C)\otimes H^*(BT)$, this shows the claim in the first case.

    In the other case, set $x=x_1, y=y_2$ and use Poincar\'e-duality to find elements $x_2, y_1$ with $x_i \cup y_i=[F]$. Choose again equivariant lifts as before. Then one has $\lr{x_i} \cup \lr{y_i}=\lr{[F]}$, in particular all elements $\lr{x_i} ,\lr{y_i}$ are nonzero. On the other hand $\lr{[F]} (\lr{y_2})_C=\lr{x_1}\lr{y_1}\lr{y_2}=0\lr{y_1} =0$ since $x_1 \cup y_2=0$. Thus $\lr{[F]}$ is a zero divisor in $H^*_T(C)$. The same equation shows the claim for $\lr{y}$, and symmetry in $x,y$ shows it for $\lr{x}$.
\end{proof}

\begin{remark}
    The following are equivalent for classes $x\in H^*_T(C)$, when the action of $T$ on $C$ is trivial:
    \begin{itemize}
        \item $x$ is a zero divisor or $x=0$.
        \item $x$ has trivial component in $H^0(C)\otimes H^*(BT)\subset H^*(C)\otimes H^*(BT)=H^*_T(C)$
        \item $x$ is not a unit in the localized ring $k\otimes_{H^*(BT)} H^*_T(C)$
        \item For any, equivalently every, point $p\in C$, the restriction $x|_p\in H^*_T(\{p\})$ is zero.
    \end{itemize}

    The last condition is the most geometric, but we will usually work in the first two formulations since they are easier to work with directly.
\end{remark}

\begin{lemma}\label{nounits}
    Assume $C$ is such that $\lr{[F]}$ is a zero divisor. Then for any class $x\in H^i(F)$ with $i$ dividing $\dim(F)$, $\lr{x}$ is a zero divisor.
\end{lemma}
\begin{proof}
    By assumption on $i$, $x^k=c[F]$ for some $k$ and (possibly zero) constant $c$. Let $0\ne y_C\in H^*_T(C)$ be such that $\lr{[F]}\cup y_C=0$. Then $(\lr{x})^ky_C=c0=0$, so $\lr{x}$ is a zero divisor. Note that $\lr{x}$ cannot vanish since for any $y$ Poincar\'e-dual to $x$, the equation $xy=[F]$ gives $\lr{x}\lr{y}=\lr{[F]}\ne 0$.
\end{proof}

The preceding lemmas show that it is quite easy for cohomology classes of $F$ to end up as zero divisors in $H^*_T(C)$. The main curvature input is the following, showing that the class $p_1(\Nu_F)$ must somehow dodge all these restrictions:

\begin{lemma}\label{frankelcontradiction}
    $\lr{p_1(\Nu_F)}$ is not a zero divisor for any fixed point component $C\ne F$.
\end{lemma}

\begin{proof}
    Assume it is. By Lemma \ref{p1lemma}, the equivariant lift of $p_1(\Nu_F)\otimes 1$ is given by $p_1^T(\Nu)-\sum \lift{x_{k\beta}}\otimes k\beta - 1\otimes C\beta^2$, where $\lift{x_{k\beta}}$ are equivariant lifts of some classes $x_{k\beta}\otimes 1$ corresponding to weights $k\beta$, produced by the middle terms in Lemma \ref{p1lemma}, and $p_1^T(\Nu)$ is the genuine equivariant first Chern class of the global bundle $\Nu$ on $F_\alpha$. Now by assumption, this class has trivial restriction to $H^0(C)\otimes H^4(BT)$. Direct computation, using Lemma \ref{nounits} to see that $\lr{x_{k\beta}}$ has no part in $H^0(C)\otimes H^2(BT)$ and Lemma \ref{p1lemma} to compute $p_1^T(\Nu)$ restricted to $C$, gives 
    $$0=\left (p_1^T(\Nu)-\sum \lift{x_{k\beta}}\otimes k\beta - 1\otimes C\beta^2\right )_{H^0\otimes H^4}=\sum_{\chi_C} 1\otimes \chi^2_C - 1\otimes C\beta^2,$$
    where the $\chi_C$ range over the weights of $T^2$ on $\Nu|_C$. Thus the sum of squares $\sum \chi_C^2$ is a multiple of $\beta^2$. This can only happen if each $\chi_C$ is itself a multiple of $\beta$ (write the sum in any basis extending $\beta$ and notice that other quadratic terms cannot cancel).

    Consider now the fixed point component $C_\beta$ of the action of $S^1_\beta$ on $M'$ containing $C$, and the corresponding component $F_\beta$ containing $F$. Since $S^1_\beta$ acts with the weights of $T^2$ quotiented out by setting $\beta=0$, the tangent space of $C_\beta$ at $C$ has dimension $\dim(C)+2b$. The dimension of $F_\beta$ is $\dim(F)+2b$. Now, on the one hand, 
    $$\dim(F_\beta)+\dim(C_\beta)=\dim(F)+\dim(C)+4b\ge\dim(F)+2b+2a=\dim(M'),$$
    so $F_\beta$ and $C_\beta$ intersect by Frankel. On the other hand, fixed point components intersect only if they are equal, and $C_\beta\cap F_\alpha \supset C \nsubseteq F_\beta\cap F_\alpha=F$. This is a contradiction.
\end{proof}

\begin{cor}\label{cor:4per}
    $H^*(F)$ is $4$-periodic with $4$-periodicity induced by $p_1(\Nu_F)$.
\end{cor}

\begin{proof}
    if $F$ is a cohomology sphere it is $4$-periodic vacuously.
    Now to establish $4$-periodicity it is sufficient to prove that $\cup p_1(\Nu_F)$ is injective in all degrees $0\le *\le \dim(F)-4$. Then Poincar\'e-duality forces the dimensions of $H^*$ and $H^{*+4}$ to be equal, which upgrades injectivity to isomorphism.
    To see this, first note that $p_1(\Nu_F)=0$ would be in contradiction to Lemma \ref{frankelcontradiction}. Now if $p_1(\Nu_F) \cup y=0$ had any nontrivial solution $y$, the last claim of Lemma \ref{Fzerodivisor} would also put us in contradiction to Lemma \ref{frankelcontradiction}. 
\end{proof}

\begin{cor}\label{cor:zerodivclass}
    The first situation of Lemma \ref{Fzerodivisor} cannot occur. The second can occur only if $F$ has the rational cohomology of $S^2\times \hp^k$.
\end{cor}

\begin{proof}
    Lemma \ref{nounits} and \ref{frankelcontradiction} show $\dim(F)=2$ mod $4$. Thus by the previous Corollary $H^*(F)$ is $4$-periodic and we may factor $[F]=p_1(\Nu_F)^k \cup x$, where $x\in H^2(F)$ is Poincar\'e-dual to $p_1(\Nu_F)^k$. Moving to $C$, we get $(\lr{p_1(\Nu_F)})^k\cup \lr{x}=\lr{[F]}$. If now we are in the first case of Lemma \ref{Fzerodivisor}, the right hand side lives in $H^{\ge 3}(C)\otimes H^*(BT)$. If $\lr{p_1(\Nu_F)}$ had any part in $H^0(C)\otimes H^4(BT)$, the left side would see $\lr{x}$ in $H^{\le2}(C)\otimes H^*(BT)$, so either $\lr{x}=0$ or $\lr{(p_1(\Nu_F)}$ is a zero divisor. Both are contradictions. Thus $H^*(F)$ is $4$-periodic with no odd cohomology. The only $4$-periodic cohomology ring with trivial odd Betti numbers having classes that multiply to zero is that of $S^2\times \hp^k$, cf. \cite[Proposition 1.5]{Kennard2021-ge}.
\end{proof}

\begin{cor}
    Theorem \ref{main:t2} holds.
\end{cor}

\begin{proof}
    The non-spherical $4$-periodic cohomology rings where $xy=0$ has no solutions are the cohomology rings of $\cp^k$ and $\hp^k$. (This follows from Poincar\'e-duality and the fact that any divisor of a class inducing periodicty induces periodicity itself, see again \cite[Proposition 1.5]{Kennard2021-ge}). Together with the possibility that $F$ is a cohomology sphere, this gives exactly the allowed outcomes of Theorem \ref{main:Fixedpointclassification} that are not covered by Lemma \ref{Fzerodivisor} and Corollary \ref{cor:zerodivclass}. 
\end{proof}

\begin{cor}
    Theorem \ref{main:Fixedpointclassification} holds.
\end{cor}

\begin{proof}
    Theorem \ref{main:t2} $\implies$ Theorem \ref{main:Fixedpointclassification} is the content of Reduction $\Delta$.
\end{proof}

\begin{cor}
    Theorem \ref{main:Hopf} holds.
\end{cor}
\begin{proof}
    Let $\tilde{M}$ be the orientation cover of $M$. Since $\chi(\tilde{M})=2\chi(M)$, we may assume $M$ is orientable. It is a result of Lefschetz that the Euler characteristic of $M$ agrees with that of the fixed point set $M^T$. This fixed point set cannot be empty due to a result of Berger. Finally, by Theorem \ref{main:Fixedpointclassification}, every fixed point component has positive contribution.
\end{proof}

\section{Orbifolds}\label{sec:orbifolds}
Generalizing the results to orbifolds requires a notion of ineffective orbifolds. To avoid technicalities arising from the groupoid perspective, we use Pardon's result that any orbifold is a global quotient as the definition \cite{Pardon_2022} and refer to \cite{MoerdijkPronk97} for the groupoid discussion.\\
An orbifold $\OO$ is defined to be an almost free action of a compact Lie group $G$ on a manifold $M$. We denote $\OO = M/G$ for the orbifold and $\vert \OO\vert$ for the topological quotient, i.e. the underlying space of $\OO$. $\OO$ is called effective if the principal isotropy group is trivial and ineffective otherwise. If the orbifold is effective, then this is equivalent to the definition in charts: Choose a Riemannian metric on $\OO$, then a classical result states that the orthogonal frame bundle $\Fr\OO$ is a manifold, such that $\O(n)$ acts almost freely on $\OO$ and that $\OO \cong \Fr\OO/\O(n)$. \\
An orbi vector bundle $E \to \OO$ over a possibly ineffective orbifold $\OO = M/G$ is given by a $G$-equivariant vector bundle $\tilde E \to M$, such that $E \cong \tilde E/G$. If $E$ (but not necessarily $\OO$) is effective, then (after choosing a suitable Riemannian metric) this induces a vector bundle $\Fr E \to \Fr E\vert_{\OO}$, such that $\OO \cong \Fr E\vert _\OO /\O(n)$.

For $\OO = M/G$, the topological orbifold-invariants of $\OO$ are defined to be the invariants of the Borel-construction $B \OO \coloneqq M_G = EG \times _G M$. If $\OO$ is effective, then $BO = \Fr\OO_{\O(n)}$. If $\OO$ is effective and  a Lie group $L$ acts on $\OO$, the induced action on $\Fr\OO$ induces an action on $B\OO$. Hence we define the $L$-equivariant topological invariants of $L \acton\OO$ to be the ones of $B\OO_L$. \\
For an orbivector bundle $E \to \OO  = M/G$, the bundle $\tilde E \to M$ induces a topological vector bundle $\tilde E_G \to M_G  = B\OO$. Hence the orbifold characteristic classes are defined to be the ones of this vector bundle. 

In order to generalize the proofs of the main theorems for orbifolds, we need orbifold analogues of the following:
\begin{multicols}{2}
    \begin{enumerate}
    \item Rational Poincar\'e duality,
    \item\label{fixcomp} $\chi(\vert\OO\vert) = \chi(\vert \OO \vert^T)$  and \\$\chi(\vert \OO^\orient\vert) = 2 \chi(\vert \OO^{2m}\vert)$,
    \item\label{conn} The Connectedness Lemma,
    \item\label{p1c1lemmas} Lemmas \ref{p1lemma} and \ref{orientationlifting}\\ \mbox{}
\end{enumerate}
\end{multicols}

It is a classical result that $H_\orb^\ast (\OO,\Q) = H^\ast(\vert \OO\vert, \Q)$ fulfills Poincar\'e duality. Furthermore, $\chi(\vert \OO^\orient\vert) = 2\chi(\OO)$, where $\OO^\orient$ denotes the orientation cover easily follows from poincare duality and the fact that $\OO = \OO^\orient/\langle \tau \rangle$ for an orientation reversing deck transformation $\tau$.

The second statement follows from \cite[Theorem~1.1]{SIKORA2004725}. For Lemma \ref{orientationlifting} note that there exists an Atiyah-Bott localization theorem for orbifolds after including an invertible factor in each term \cite{MEINRENKEN1998240}. Hence the same proof as in Lemma \ref{orientationlifting} applies. For Lemma \ref{p1lemma} the same proof applies, by the following 
\begin{lemma}\label{lem:trivialaction}
    Let $E \to F$ be an  Riemannian orbi vector bundle with $E$ effective. Let $T\acton E$ fibrewise, i.e. fixing $F$. Then $$H_{T,\orb}^\ast(F, \bb Q) \cong H^\ast(\vert F\vert,\bb Q) \tensor H^\ast (BT, \bb Q)$$ 
\end{lemma}

\begin{proof}
   There is a natural $T$-action on $\Fr E$, coming from differentiating the $T$-action on $E$. It is well known that the projection $q \colon BF = (\Fr E\vert _F)_\O(n) \to  \vert F\vert$ gives an isomorphism of cohomology rings $q^\ast \colon H^\ast(\vert F\vert , \bb Q)\to  H^\ast_{orb} (F,\bb Q)$. Since $q$ is equivariant, it induces a map $q_T \colon BF_T \to \vert F\vert_T = BT \times \vert F \vert$, which is a map of fibrations
  \begin{center}
\begin{tikzcd}
    BF \arrow[r] \arrow[d, "q"] & BF_T \arrow[r] \arrow[d,"q_T"] & B T \arrow[d, "="] \\
    \vert F \vert  \arrow[r]   & B T \times \vert F \vert  \arrow[r]         & BT.
\end{tikzcd}
    \end{center}
   Since $BT$ is simply connected, $q $ induces an isomorphism on the $E_2$ page of the Serre spectral sequence with rational coefficients
   $$H^p(BT,\bb Q) \tensor H^q(\vert F \vert,\bb Q) \to H^p(BT, \bb Q) \tensor H^q(BF_T, \bb Q), x \tensor y \mapsto x \tensor q^\ast y.$$
    By Zeeman's comparison theorem for spectral sequences \cite{Zeeman1957-ik}, $q^\ast_T$ is also an isomorphism, which proves the claim.
\end{proof}

Problem \ref{conn} is solved in the following generalization of the Connectedness Lemma:

\begin{Lemma} \label{lem:connectedness}
    Let $(\OO^n,g)$ an effective Riemannian orbifold with positive sectional curvature.
    \begin{enumerate}
        \item  $N^{n-k} \subset \OO$ a compact totally geodesic embedded suborbifold with orbifold structure $N =\Fr\OO\vert_N /\O(n)$. Then the inclusion $N\into \OO$ is orbifold $(n-2k+1)$-connected. 
        \item Let $N^{n-k_1}_1$ and $N_2^{n-k_2}$ be compact totally geodesic embedded suborbifolds, such that $k_1 \le k_2$, $k_1 + k_2 \le n$. Then the intersection $N_1^{n-k_1} \cap N_2^{n-k_2}$ is a totally geodesic embedded suborbifold as well and the inclusion
        \begin{align*}
            N_1^{n-k_1} \cap N_2^{n-k_2} \to N_2^{n-k_2}
        \end{align*}
        is orbifold $(n-k_1-k_2)$-connected.
    \end{enumerate}
\end{Lemma}

\begin{remark}
    Even though $\OO$ is assumed effective, the same follows for neither $N$ nor any $N_i$, since the action of $O(n)$ on $Fr \OO|_N$ may have principal isotropy. In particular the conclusions of the lemma are not guaranteed to hold when suborbifolds are defined as the effectivization of the induced orbifold structure.
\end{remark}

For the proof, we need the following:
\begin{Lemma}\label{lem:Framemetric}
    There is a metric on $\Fr\OO$ such that the quotient map $\pi \colon \Fr \OO \to \OO$ is a Riemannian submersion with parallel fibres and all horizontal sectional curvatures are positive.
\end{Lemma}
\begin{proof}
   We denote by $\omega \in \Omega(\Fr\OO) \tensor \mathfrak o (n)$ the Riemannian connection $1$-form of $(\OO,g)$. $\V = \ker d \pi$ denotes the vertical and $\H = \ker w$ the horizontal distribution. Note that $\V_f \cong \mathfrak o(n)$ for a frame $f \in \Fr \OO$.\\
   We define the following metric on $\Fr \OO$: For $v,w \in \mathfrak{o}(n)$  set $Q(v,w) = -\Tr(v\cdot w)$. Let $a>0$. For $v,w \in \V$ and $x,y \in \H$ we set
   $$g_\Fr(v,w) = a^2 Q(\omega(v), \omega(w)), \quad g_\Fr(v,x) =0, \quad g_\Fr(x,y) = g(d\pi(x), d\pi(y))$$
   With this definition the fibres of $\pi$ are parallel. Now at a frame $f \in \Fr(\OO)$ we obtain from O'Neills result (see \cite[section 5]{ONeill1966}):
   \begin{align*}
    \sec(x\wedge y)_{g_\Fr} &= \sec_g(d\pi(x)\wedge d\pi(y)) - \frac{3 a^2\cdot\norm{R^g(d\pi(x)\wedge d \pi(y))}^2}{4\norm{x \wedge y}^2}\label{sec_horizontal}
\end{align*}
where $A$ denotes the O'Neill-tensor and $\alpha = \sum \omega_i^j(v)f_i\wedge f_j$ locally. Note that this argument works for orbifolds, since you can apply it to a local trivialisation of $\Fr\OO$. Since the sectional curvature of $\OO$ is positive, we obtain that for sufficiently small $a$, the sectional curvature $\sec_{g_\Fr}(x\wedge y)$ is positive. 
 \end{proof}

\begin{proof}[Proof of Lemma \ref{lem:connectedness}] We follow the proof of \cite{Wil03}, but lift the arguments to the frame bundle ($\Fr\OO$, $g_\Fr$) with the metric obtained from Lemma \ref{lem:Framemetric}, which has positive horizontal sectional curvature.
 We denote $\tilde N = \Fr \OO\vert_N$. We consider the path space $\Omega_{\tilde N}(\Fr\OO)$ of piecewise smooth paths in $\Fr\OO$ starting and ending in $\tilde N$ and give a lower bound for the index of critical points of the energy functional $E$. A critical point of $E$ is a geodesic starting and ending perpendicularly to $\tilde N$. Let $c$ be such a geodesic. Note that $c$ is always perpendicular to the $\O(n)$-orbits.
    Consider the following family $V = \lbrace X \setvert X(t) \in  \H, X^\prime(t) \in \V\rbrace$ of vector fields.
    Obviously $\dim V = n$, since these are horizontal lifts of parallel fields in $\OO$. Now consider the subfamily
    $W = \lbrace X \in V\setvert X \perp c,   X(0) \in T_{c(0)} \tilde N, X(1) \in T_{c(1)}\tilde N, X^\prime(0) = 0\rbrace$, which clearly has
     $\dim W \ge n -2k  + 1$. 
    Therefore we have an at least $(n-2k+1)$-dimensional family of horizontal vector fields, which are normal to $c$, starting and ending tangentially to $\tilde N$. Since $\tilde N$ is horizontally totally geodesic, by choice of a suitable $a$ for $g_\Fr$ the second variation of energy implies that the index of $c$ is at least $n-2k+1$ by the positivity of the  horizontal sectional curvature of $\Fr\OO$. By following the argument of \cite{Wil03}, we obtain that the inclusion $\tilde N \to \Fr\OO$ is $(n -2k +1)$-connected. 
    The inclusion induces a map of principal bundles $(\O(n) \to E\O(n)\times \Fr\OO \to B \OO) \to (\O(n) \to E\O(n)\times \tilde N \to B N)$
    This induces a map between the long exact sequences of homotopy groups. Since $E\O(n)$ is contractible, the result follows from the five lemma. With this setup the second statement follows in the same way as in \cite{Wil03}.
\end{proof}
\begin{remark}
    The Connectedness Lemma is in generally false for the underlying spaces: Let $N = \bb S^{2n+1}/\bb Z_m$ a lense space and $\OO =  \bb S^{2n+2}/\bb Z_m \cong \Sigma(N)$ its suspension. Then $N$ is a totally geodesic sub-orbifold of $\OO$ of codimension 1, but $ \vert N \vert \to \vert \OO \vert$ is not even $2$-connected, since $\pi_1(\vert N \vert) = \bb Z_m$ and $\pi_1(\OO) = 0$.
\end{remark}

\begin{appendices}
 
\section{An approach via the first
Chern class}\label{sec:c1}

In this appendix, we want to describe the originally found way to force odd cohomology to vanish, mostly because it controls a different quantity (the sum of the weights instead of the sum of the squares of the weights) and works in slightly different situations, making it very likely both approaches will be useful in future research. The main difference in idea, apart from using a different characteristic class, is that one proceeds by induction:\\
\noindent
\textbf{Reduction} $\square$: Assume $M$ has minimal dimension among even-dimensional positively curved spaces with $T^3$-action containing a fixed point component with nontrivial odd rational cohomology.

\begin{lemma}
    In the situation of Reduction $\square$, the action of $T^3$ has no finite isotropy groups near any fixed point component with nontrivial odd cohomology.
\end{lemma}

\begin{proof}
    If $I$ is a finite isotropy group near some $F$, then $\tilde T^3=T^3/I \acton Fix(I)$ is an example of smaller dimension that still contains $F$.
\end{proof}

After this, pick again a component $F$ and perform Reduction $\Delta$ as before to replace the $T^3$-action on $M$ by a $T^2$-action that splits at $F$. Now instead of all weights at $F$ having the form $k\beta$, $l\alpha$ for potentially varying $k,l$, the isotropy condition forces $k=l=1$. Indeed, if any weight is divisible by an integer $\ge 2$, one can find a finite isotropy group near $F$ whose order is also divisible by that integer. Since weights of real representations are defined only up to sign, this shows the classes of $\beta$ and $\alpha$ are the only that occur. In fact, in the usual process of using the torus action to define a complex structure on the bundle, the obtained complex weights all end up in a chosen open half-space, and that half space may be chosen to be one containing $\alpha, \beta$ to get only positive signs. Let $a, b$ denote the multiplicities of $\alpha$ resp. $\beta$.

Now the argument goes like this:

\begin{lemma}
    If $F$ is simply connected with nontrivial odd cohomology, then the equation $c_1(\Nu|_F) \cup y=0$ has a solution for some $y\ne0$ with odd degree.
\end{lemma}

\begin{proof}
    Being simply connected forces $H^1(F)=H^{\dim F-1}(F)=0$. Thus taking any odd class, repeated product with $c_1\in H^2(F)$ eventually gives zero. Take $y$ to be the last nonzero class obtained along this process.
\end{proof}

We may again shift this to a second component $C$ using a simplified version of Lemma \ref{Fzerodivisor}:

\begin{lemma}
    There is another component $C$ where $\lr{y}$, the restriction of the equivariant lift of $y\otimes 1$, is nonzero.
\end{lemma}
\begin{proof}
    Take a Poincar\'e-dual $z$ of $y$, and consider $\lift{y} \cup\lift{z}$. Atiyah-Bott localization gives a nonzero contribution at $F$ and a global value of zero for the integral, so there is a second component with nonzero contribution. In particular, $\lr{y}$ and $\lr{z}$ cannot vanish there. 
\end{proof}

Finally, we look at the equation $c_1(\nu|_F)y=0$ at the component $C$. For this, first we need to know what the equivariant Chern class looks like:

\begin{lemma}\label{c1lemma}
    Let $V\to F$ be a $T$-equivariant complex vector bundle with the $T$-action fixing $F$ pointwise. Then the equivariant first Chern class satisfies $$c_1^T(V)=c_1(V)\otimes 1 + \sum 1 \otimes \chi\in H^*(F)\otimes H^*(BT)=H_T^*(F),$$ where $c_1(V)$ is the non-equivariant first Chern class and the sum ranges over the weights of all irreducible subrepresentations of $V$ at any point $p\in F$ (with multiplicity), with the lattice of weights identified with $H^2(BT, \Z)\subset H^2(BT, \Q)$.
\end{lemma}

\begin{proof}
    Use the weights to split $V$ into a direct sum of subbundles $V_\chi$ on which $T$ acts by $\chi$ at every point. Denote by $V_\chi^0$ the same bundles but with $T$ acting trivially. Then $T$-equivariantly we have $V_\chi=V_\chi^0\otimes_\C L_\chi$, where $L_\chi=F\times \C \to F$ is the trivial line bundle with $T$ acting on $\C$ via $\chi$. Since $L_\chi$ is a line bundle, $c_1^T(V_\chi)=c_1^T(V_\chi^0)+\dim_\C(V_\chi^0)c_1^T(L_\chi)$. The bundle $L_\chi$ has $c_1^T(L_\chi)=1\otimes\chi$, thus by additivity of the first Chern class we get $c_1^T(V)=\sum_\chi c_1^T(V_\chi^0) + \sum_\chi \dim_\C(V_\chi)\chi =  c_1(V)\otimes 1 + \sum_\chi \dim_\C(V_\chi)\chi$, which is the claimed formula, noting that $c_1^T=c_1\otimes 1$ for bundles with trivial actions.
\end{proof}

This lets us relate the representation data at $C$ with that at $F$:

\begin{lemma}
    At $C$, we have $\sum\chi_i^C=b\beta$.
\end{lemma}
\begin{proof}
    The equation $c_1(\nu|_F)y=0$ holds at $F$ and lives in degree $\le \dim(F)-1$, so in the range where $H^*_T(F_\alpha)\to H^*_T(F)$ is an isomorphism. Thus the equation remains true after lifting both classes to $H^*_T(F_\alpha)$, and then after restricting to $C$, and we get $\lr{c_1(\Nu|_F)} \cup \lr{y}=0$. In particular, since $\lr{y}\ne 0$, $\lr{c_1(\Nu|_F)}$ is a zero divisor in $H^*_T(C)$. This means the component of this class in $H^0(C)\otimes H^2(BT)$ must vanish. Now we can use Lemma \ref{c1lemma} to compute this component, and get $b\beta-\sum\chi_i^C=0$, where the sum ranges over all weights of the action of $T^2$ on $\Nu|_C$ (and the other contribution is that same sum for $\Nu|_F$).
\end{proof}

Thus the first Chern class controls the sum of the weights, where the first Pontryagin class controlled the sum of the squares. However, it turns out with the finite isotropy reduction, this is just as good:

\begin{lemma}\label{weightsatC}
    The $T^2$-action on $\Nu|_C$ has $\beta$ as the only weight, occurring with multiplicity $b$. 
\end{lemma}

\begin{proof}
    Let $\chi_i$ denote the weights of $T^2$ on $\Nu|_C$. 

    Consider the restricted action of the circle $S^1_\alpha\subset T^2$. It acts via the weights of the $T^2$-action after quotienting out by $\alpha=0$. All of $F_\alpha$ is a single $S^1_\alpha$-fixed point component, so the weights of $S^1_\alpha$ on the normal bundle $\Nu$ to $F_\alpha$ are locally constant. Thus we get that the $T^2$-weights at $C$ agree with those at $F$ after quotienting out by $\alpha$, so $\chi_i=\beta+k_i\alpha$ for some integers $k_i$.

    At this point we make use that we have Reduction $\square$ also at the component $C$, since $0\ne \lr{y}\in H^{odd}_T(C)=H^{odd}(C)\otimes H^*(BT)$ means $C$ has odd cohomology. Consider two different weights $\chi_i, \chi_j$. They span some sublattice $\Z[\chi_i, \chi_j]$ of the weight lattice of $T$, and the common kernel of this sublattice is a finite subgroup $I$ of $T^2$ whose fixed point component at $p$ is strictly larger than $F$, since its tangent space at $p$ contains all the spaces with weights in $\Z[\chi_i, \chi_j]$. The points in this fixed point space have generic isotropy $I$, which means that we are in contradiction with Reduction $\square$ unless $|I|=1$. On the other hand, $|I|$ equals the index of $\Z[\chi_i, \chi_j]$ in the weight lattice, which in an integral basis of that lattice is simply the (absolute value of the) determinant of the matrix $(\chi_i ,\chi_j)$. Thus we have in basis $(\beta, \alpha)$ that $1=|\det\begin{pmatrix}
        1 & 1 \\k_i & k_j
    \end{pmatrix}|$, so $|k_i-k_j|\le 1$ for all $i,j$.

    To conclude, we know that we have $b\beta=\sum\limits_{i=1}^b\chi_i=\sum\limits_{i=1}^b\beta+k_i\alpha=b\beta+(\sum k_i)\alpha$, so $\sum k_i=0$. However, the restriction $|k_i-k_j|\le 1$ shows that no two $k_i$ can have different signs, so they must all vanish. Thus $\chi_i=\beta$ for all $i$, which was the claim.
\end{proof}

This, finally, is in direct contradiction to Frankel's theorem, by the argument we already saw in the proof of Lemma \ref{frankelcontradiction}:

\begin{cor}
    The situation of Reduction $\square$ cannot occur. In particular, all fixed point components have vanihsing odd cohomology.
\end{cor}

\begin{proof}
    Consider the fixed point component $C_\beta$ of the action of $S^1_\beta$ containing $C$, and the corresponding component $F_\beta$ containing $F$. Since $S^1_\beta$ acts with the weights of $T^2$ quotiented out by setting $\beta=0$, the tangent space of $C_\beta$ at $C$ has dimension $\dim(C)+2b$ by Lemma \ref{weightsatC}. The dimension of $F_\beta$ is $\dim(F)+2b$. Now, on the one hand, $\dim(F_\beta)+\dim(C_\beta)=\dim(F)+\dim(C)+4b\ge\dim(F)+2b+2a=\dim(M')$, so $F_\beta$ and $C_\beta$ intersect by Frankel. On the other hand, fixed point components intersect only if they are equal, and $C_\beta\cap F_\alpha \supset C \nsubseteq F_\beta\cap F_\alpha=F$. This is a contradiction.
\end{proof}

\end{appendices}

\printbibliography

\end{document}